\documentclass{article}

\usepackage{graphicx}

\title{Pencils on double coverings of curves.}
\author{Luis Fuentes Garc\'{\i}a}
\date{}
\newtheorem{teo}{Theorem}[section]

\newtheorem{cor}[teo]{Corollary}
\newtheorem{lemma}[teo]{Lemma}

\newtheorem{que}[teo]{Question}

\font\euf=eufm10 at 12pt

\def\g2{\pi}

\def\e{\mbox{\euf e}}

\def\b{\mbox{\euf b}}

\def\d{\mbox{\euf d}}

\def\c{\mbox{\euf c}}

\def\a{\mbox{\euf a}}

\def\E{{\cal E}}

\def\aa{\mbox{\euf a}}

\def\P{{\bf P}}

\newcommand\Te{{\cal O}}

\def\ZZ{\leavevmode\hbox{$\rm Z$}}

\def\impp{{\quad\Rightarrow\quad}}

\def\qed{\hspace{\fill}$\rule{2mm}{2mm}$}
\newcommand\lrw{\longrightarrow}

\begin{document}

\maketitle

{\bf Abstract:} Let $X$ be a smooth curve of genus $g$. When $\pi\geq 3g$ and $d\geq \pi-2g+1$ we show the existence of a double covering $\gamma:C\lrw X$ where $C$ a smooth curve of genus $\pi$ with a base-point-free pencil of degree $d$ which is not the pull-back of a pencil on $X$. \\


\vspace{0.1cm}

\section{Introduction.}

Let $C$ and $X$ be two smooth curves of genus respectively $\pi$ and $g$. Let $\gamma:C\lrw X$ be a double covering. It is clear that the base-point-free pencils $g_d^1$ on $X$ induce base-point-free linear systems on $C$ which have low degree respect to the genus $\pi$. However, the existence of base-point-free pencils of given degree $d$ on $C$ which are not induced by $\gamma$ is a more subtle problem.

The following facts are known:

\begin{enumerate}

\item If $d\leq \pi-2g$, a pencil on $C$ of degree $d$ is always induced by a pencil on the curve $X$. This is a consequence of the Castelnuovo-Severi inequality (see \cite{Ka84}).

\item If $X$ is a generic curve and $\pi\geq 5g+1$ then there is a base-point-free pencil of any degree $d\geq \pi-g$ which is not composed with the given double covering (see \cite{BaKe95}).

\item If $0\leq s\leq 2g$ and $\pi\geq 4g-s$ then there is a base-point-free pencil of  degree $d=\pi-2g+1+s$ which is not composed with the given double covering (see \cite{KeOh06}). Note, that these bounds implies $d\geq 2g+1$.

\item Given a smooth curve $X$ of genus $g$, if $\pi\geq 4g+5$ and $d\geq \pi-2g+1$ then there exist a double covering $\gamma:C\lrw X$ with $C$ smooth of genus $\pi$ and a base-point-free linear system of degree $d$ on $C$ which is not induced by $\gamma$ (see \cite{BaKePa04}). In this case the double covering is not arbitrary. 

\end{enumerate}

In this paper we improve significantly the bound on the genus $\pi$ for the existence of the double covering and the base-point-free pencil of degree $d$. The main result is the following:

\begin{teo}

Let $X$ be a smooth curve of genus $g$. If $\pi\geq 3g$ and $d\geq \pi-2g+1$, then there is a double covering $\gamma:C\lrw X$ with $C$ a smooth curve of genus $\pi$ and a base-point-free divisor $\d\in Div^d(X)$ such that the linear system $|\d|$ is not induced by a linear system on $X$.

\end{teo}

The proof uses the relation between double covers and $2$-secant curves on ruled surfaces. Note, that this fact was used by Keem and Ohbuchi in \cite{KeOh06}, but not exactly in the same way. They use the properties of the elementary transformations of ruled surfaces. Here, given a smooth curve $X$ and a decomposable ruled surface $S=\P(\Te_X\oplus \Te_X(\e))$, we describe with detail the construction of a $2$-secant curve $C$ on $S$ with prescribed branch points. The pencils on $C$ will be the restriction of suitable pencils of $S$.

Finally, we give some results relating the gonality and the Clifford index of the curves $X$ and $C$.

We refer to \cite{FuPe05} and \cite{Ha77} for a systematic development of the projective theory of ruled surfaces and to \cite{FuPe05*} for the study of $2$-secant curves and double coverings.

\section{Double coverings and ruled surfaces.}

A {\it geometrically ruled surface}, or simply a {\it
ruled surface}, will be a $\P^1$-bundle over a smooth curve $X$ of genus $g>0$. It will be
denoted by $\pi: S=\P(\E_0)\lrw X$ with fibre $f$. We will suppose that
$\E_0$ is a normalized sheaf and $X_0$ is the section of minimum self-intersection
that corresponds to the surjection $\E_0\lrw \Te_X(\e)\lrw 0$, with $\bigwedge^2\E\cong
\Te_X(\e)$ and  $e=-deg(\e)$. We know that  $Num(S)=\ZZ X_0\oplus \ZZ f$ (see \cite{Ha77},V.2 and \cite{FuPe05}).

When $\E_0\equiv \Te_X\oplus \Te_X(\e)$ we say that the ruled surface is decomposable. In this case we will denote by $X_1$ the irreducible curves of the linear system $|X_0-\e f|$. They does not meet $X_0$ and have self-intersection $e$.

A curve $C\subset S$ will be said $n$-secant when $C\cdot f=n$. The following results describe the relation between $2$-secant curves on decomposable ruled surfaces and double coverigns.

\begin{teo}

Let $\gamma:C\lrw X$ be a double covering of a smooth curve $X$. Let ${\cal L}$ be an invertible sheaf of $C$. Then:

\begin{enumerate}

\item $\E=\gamma_*{\cal L}$ is a local free sheaf on $X$ of rank $2$. From this, $p:\P(\E)\lrw X$ is a geometrically ruled surface.

\item There is a closed immersion $j:C\lrw \P(\E)$ verifying $p\circ j=\gamma$.

\item $H^i(\E)=H^i(\Te_{P(\E)}(1))\cong H^i({\cal L})$.

\item $deg({\cal L})=deg(\Te_{P(\E)}(1))+\pi-2g+1$.

\end{enumerate}

\end{teo}
{\bf Proof:} See \cite{FuPe05*}, 2. \qed

\begin{lemma}\label{construction}

Let $p:S=\P(\Te_X\oplus \Te_X{(\e)})\lrw X$ be a decomposable ruled surface over a smooth curve $X$. Let $\b_1,\b_2$ be two effective divisors of $X$ verifying:

\begin{enumerate}

\item $\b_1$ and $\b_2$ are smooth.

\item $\b_1\cap \b_2=\emptyset$.

\item $\b_2\sim \b_1+2\e$.

\item If $\e\sim 0$, then $\b_1\not\sim 0$.

\end{enumerate}

Then, there is a smooth curve $C\sim 2X_0+\b_1 f$ with an induced double covering $\gamma:C\lrw X$ with branch divisor $\b_1+\b_2$.

\end{lemma}
{\bf Proof:} The idea is the same of the proof of the Theorem 2.9 in \cite{FuPe05*}. Let us consider the nontrivial involution $\sigma$ of $S$ with fixed divisors $X_0$ and $X_1$. The divisors of the following pencil are invariant by the involution $\sigma$:
$$
\langle 2X_0+\b_1 f, 2X_1+\b_2 f\rangle \subset |2X_0+\b_1 f|.
$$
Let $C$ be a generic divisor of the pencil. By Bertini's Theorem, $C$ is smooth away from the base locus. Since $\b_1\cap \b_2=\emptyset$, the base locus of the pencil is $\b_1 f \cap X_1$ and $\b_2 f\cap X_0$. Then, the smoothness of $\b_1$ and $\b_2$ implies the smoothness of $C$.

Now, let us see that $C$ is irreducible. $C$ can not contain a fiber $f$, because the points of $(C-f)\cap f$ would be singular points. On the other hand, suppose that $C=A+B$, where $A$ and $B$ are irreducible $1$-secant curves. Since $C$ is smooth, $A\cap B=\emptyset$. This implies $A\sim X_0$, $B\sim X_1$ and $\b_1\sim \b_2\sim \e\sim 0$. But this contradicts the hypothesis.

Finally, the projection $p$ induces a  double covering $\gamma:C\lrw X$. Moreover $C$ is invariant by the involution $\sigma$, so the ramification points of $\gamma$ are $C\cap X_0$ and $C\cap X_1$ and the branch divisor is $\b_1+\b_2$. \qed

\begin{lemma}\label{bound0}

Let $C$ be a smooth curve of genus $\pi$ with a double covering $\gamma:C\lrw X$ over a smooth curve $X$ of genus $g$. If $C$ is a $2$-secant curve of a ruled surface $p:S\lrw X$ with invariant $e$, then:
$$
e\leq \pi-2g+1.
$$

\end{lemma}
{\bf Proof:} Let us suppose $C\subset S$, with $C\sim 2X_0+\b f$. Because $C$ and $X_0$ are distinct curves:
$$
0\leq C\cdot X_0=deg(\b)-2e \impp deg(\b)-e\geq e.
$$
Then:
$$
\pi=g(C)=g(X_0)+g(X_0+\b f)+X_0\cdot (X_0+\b f)-1=2g+deg(\b)-e-1
$$
From this:
$$
\pi-2g+1=deg(\b)-e\geq e.
$$ \qed

\begin{teo}\label{existence}

Let $C$ be a smooth curve of genus $\pi$ with a double covering $\gamma:C\lrw X$ over a smooth curve of $C$ of genus $g$. Let $\c$ be the branch divisor. If $0\leq e\leq \pi-2g+1$ then there is a decomposable ruled surface $S$ with invariant $e$ such that $C$ is a $2$-secant curve. In particular:
$$
C\sim 2X_0+\b_1 f ,\hbox{ where }2\b_1+2\e\sim \c.
$$

\end{teo}
{\bf Proof:} It is well known that the double covering is completely determined by fixing the branch divisor. Let us choose two effective divisors $\b_1,\b_2\in Div(C)$ such that:
$$
\c=\b_1+\b_2,\quad  \hbox{ and }\quad deg(\b_1)-deg(\b_2)=2e.
$$
The map $\times 2:Pic^e(X)\lrw Pic^{2e}(X)$ is a surjection, so there is a divisor $\e$ such that $2\e\sim \b_2-\b_1$. Now, applying the  Theorem \ref{construction} we construct a double covering with the prescribed branch divisor on a decomposable ruled surface with invariant $e$.  \qed

\section{Pencils on double coverings.}

The main tool to construct double coverings with base-point-free pencils of low degree will be the following theorem. The pencils will be the restriction of unisecant linear systems on a decomposable ruled surface to $2$-secant curves. 

\begin{teo}\label{maintool}

Let $X$ be a smooth curve of genus $g$. Let $\a_1,\a_2,\b$ be effective divisors on $X$ verifying:

\begin{enumerate}

\item $\b$, $\b+2\a_2-2\a_1$ are smooth divisors without common base points.

\item $Bs|\a_1|\cap Bs|\a_2|=Bs|\a_1|\cap Bs|\b|=Bs|\a_2|\cap Bs|\b+2\a_2-2\a_1|=\emptyset$.

\end{enumerate}

Then there is a double covering $\gamma:C\lrw X$ where $C$ is a smooth curve of genus $\pi=2g-1+deg(\b+\a_2-\a_1)$ and there is a base-point-free linear system $|\d|$ in $C$ of degree $deg(\b+\a_2)$.

Moreover, if $deg(\b+2\a_2-2\a_1)>0$ then $|\d|$ is not the pull-back of a linear system in $X$.

\end{teo}
{\bf Proof:} We apply the Lemma \ref{construction} for $\e\sim \a_2-\a_1$, $\b_1\sim \b$ and $\b_2\sim \b+2\a_2-2\a_1$. We obtain a smooth curve $C$ of genus $2g-1+deg(\b+\a_2-\a_1)$ and a double cover $\gamma:C\lrw X$. Moreover, consider the linear system $|H|=|X_0+\a_1 |$ on the ruled surface $\P(\Te_X\oplus \Te_X{(\e)})$. Let us see that the restriction  $|\d|=|H|_C$ of this linear system to the curve $C$ is base-point-free.

Since $\a_1$ and $\a_2$ have not common base points, the unique base points of $|H|$ are at $X_0\cap \a_2 f$ and $X_1\cap \a_1 f$. We also know that $C\sim 2X_0+\b f$. Thus, $C\cap X_0=(\b+2\a_2-2\a_1)f\cap X_0$ and $C\cap X_1=\b  f\cap X_1$. But the pairs ($\b+2\a_2-2\a_1$, $\a_2$) and ($\b,\a_1)$ have not common base points. We deduce that $|\d|$ is base-point-free.

Finally, since $C\cdot X_0=deg(\b+2\a_2-2\a_1)>0$ the linear system $|\d|$ is not the pull-back of a linear system in $X$.  \qed

\begin{teo}\label{fundamental}

Let $X$ be a smooth curve of genus $g$. If $\pi\geq 3g$ and $d\geq \pi-2g+1$, then there is a double covering $\gamma:C\lrw X$ with $C$ a smooth curve of genus $\pi$ and a base-point-free divisor $\d\in Div^d(X)$ such that the linear system $|\d|$ is not induced by a linear system on $X$.

\end{teo}
{\bf Proof:} We will apply the Theorem \ref{maintool}. We only have to choose suitable divisors $\a_1,\a_2$ and $\b$. We consider two cases:

\begin{enumerate}

\item $d-(\pi-2g+1)$ is even.

Let $\a_1,\a_2$ be two effective divisors on $X$ of degree $a=(d-(\pi-2g+1))/2$ and without common base points. In particular, if $d=(\pi-2g+1)$ we will take $\a_1=\a_2=0$. Let $\c$ be a generic divisor on $X$ of degree $d$. By hypothesis, $d\geq g+1+2a$ so the divisors $\b=\c-2\a_2$ and $\b+2\a_2-2\a_1$ are effective and base-point-free. 

\item $d-(\pi-2g+1)$ is odd.

Let $a=(d+1-(\pi-2g+1))/2$. Let $\a_1,\a_2$ be two effective divisors on $X$ of degrees $a$ and $a-1$ respectively, without common base points. In particular, if $d=(\pi-2g+1)+1$ we will take $\a_2=0$. Let $\c$ be a generic divisor on $X$ of degree $d$. Now $d\geq g+2a$ and the divisors $\b=\c-2\a_2$ and $\b+2\a_2-2\a_1$ are effective and without common base points.

\end{enumerate}

\qed

\begin{teo}\label{lowdegree}

Let $X$ be a smooth curve of genus $g$ with a base-point-free pencil $|\b|$ of degree $d>0$. Then there is a double covering $\gamma:C\lrw X$ with $C$ a smooth curve of genus $\pi=2g-1+d$ and a base-point free-divisor $\d\in Div^d(X)$ such that the linear system $|\d|$ is not induced by a linear system on $X$.

\end{teo}
{\bf Proof:} It is sufficient to apply the Theorem \ref{maintool} for the divisors $\b$ and $\a_1=\a_2=0$ . \qed

As a consequence of this result we see that the bound $\pi\geq 3g$ of the Theorem \ref{fundamental} is not optimal. If $X$ is a curve of gonality $gon(X)=d$ we can apply the Theorem \ref{lowdegree} to obtain a $C$ with genus $\pi=2g-1+d$ and a base-point-free pencil of degree $d$. On the other hand, if $d=\pi-2g+1$ the bound $\pi\geq 2g-1+gon(X)$ can not be improved because $gon(C)\geq gon(X)$. Thus, the remaining open question is the following:

\begin{que}

Let $X$ be a smooth curve of genus $g$. Let $\pi,d$ be integers with $2g-1\leq \pi<3g$ and $d\geq max\{\pi-2g+1,gon(X)\}$. Is there any covering $\gamma:C\lrw X$ with $C$ a smooth curve of genus $\pi$ and a base-point-free pencil of degree $d$ which is not the pull-back of a pencil of $X$?.

\end{que}

\section{Gonality and Clifford index of double coverings.}

Given a double covering $\gamma:C\lrw X$, in this section we study the relation bewteen the gonality and the Clifford index of the curves $C$ and $X$. Although some of the proofs could be shorted by using the Castelnuvo-Severi inequality, we only apply the methods described in this article.

\begin{lemma}\label{bound1}

Let $\gamma:C\lrw X$ be a double covering of a smooth curve $X$. Then:
$$
gon(X)\leq gon(C)\leq 2gon(X).
$$

\end{lemma}
{\bf Proof.} Let $\d$ be a divisor of $C$ providing the gonality $gon(C)$. Then $\gamma_*\d$ is a divisor of $C$ with
$deg(\gamma_*\d)=gon(C)$ and $h^0(\Te_X(\gamma_*\d))\geq 2$. Thus, $gon(X)\leq gon(C)$.

On the other hand, by the Theorem \ref{existence} we know that there is a ruled surface $S=\P(\E)$ with invariant $e=\pi-2g+1$ containing the curve $C$. Furthermore $C\sim 2X_0+\b f$, where $\b$ is the branch divisor of $\gamma$. Let $\a$ be a divisor on $X$ providing the gonality $gon(X)$. Let us consider the linear system $|H|=|X_0+\a f|$. Then:
$$
deg(\Te_C(H))=H\cdot C=-2e+2deg(\a)+deg(\b)=2deg(\a)=2gon(X).
$$
Moreover, 
$$
h^0(\Te_C(H))=h^0(\Te_S(H))=h^0(\Te_X(\a))+h^0(\Te_X(\a+\e))\geq h^0(\Te_X(\a))=2.
$$
From this, $gon(C)\geq deg(\Te_C(H))=2gon(X)$. \qed

\begin{teo}\label{gonality}

Let $C$ be a smooth curve of genus $\pi$ with a double covering $\gamma:C\lrw X$ over a smooth curve of $C$ of genus $g$. Then:
$$
gon(C)\geq min\{\pi-2g+1,2gon(X)\}
$$

\end{teo}
{\bf Proof:}  Let $\d$ be a divisor on $C$ providing the gonality $gon(C)$. Let $(S,H)$ be the polarized ruled surface with $S=\P(\gamma_*\Te_C(\d)))$ and $\Te_S(H)=\Te_S(1)$. Let $H\sim X_0+\a f$. We know that 
$$
2=h^0(\Te_C(\d))\leq h^0(\Te_X(\a))+h^0(\Te_X(\a+\e))\qquad \hbox { and }\qquad deg(\a)\geq 0.
$$
Moreover:
$$
gon(C)=deg(\d)=H^2+\pi-2g+1=2deg(\a)-e+\pi-2g+1.
$$
Now, we distinguish two cases:

\begin{enumerate}

\item If $2deg(\a)-e\geq 0$ then: 
$$
gon(C)=2deg(\a)-e+\pi-2g+1\geq \pi-2g+1.
$$

\item If $2deg(\a)-e<0$ then $deg(\a)-e<0$. From this, $h^0(\Te_X(\a+\e))=0$ and  $h^0(\Te_X(\a))\geq 2$. We see that $deg(\a)\geq gon(X)$ and
$$
gon(C)=2deg(\a)-e+\pi-2g+1\geq 2gon(X).
$$ \qed

\end{enumerate}

\begin{teo}\label{clifford}

Let $C$ be a smooth curve of genus $\pi$ with a double covering $\gamma:C\lrw X$ over a smooth curve of $C$ of genus $g$. Then:
$$
Cliff(C)\geq min\{\pi-2g-1,2gon(X)-2\}.
$$
\end{teo}
{\bf Proof:} If the Clifford dimension of $C$ is $1$, then the bound follows from the previous theorem. Thus, let $\d$ be a divisor on $C$ providing the Clifford index $Cliff(C)$ with $h^0(\Te_C(\d))\geq 3$. We construct the polarized ruled surface $(S,H)$ corresponding to $\gamma_*\Te_C(\d)$. Note that:
$$
deg(\d)\leq \pi-1 \impp H^2+\pi-2g+1\leq \pi-1\impp H^2\leq 2g-2.
$$

We know that:
$$
\begin{array}{rl}
{Cliff(C)}&{=2deg(a)-e+\pi-2g+1-(2h^0(\Te_C(L))-2)\geq}\\
{}&{\geq 2deg(a)-e+\pi-2g-2h^0(\Te_X(\a))-2h^0(\Te_X(\a+\e))+3.}\\
\end{array}
$$
We distinguish the following cases:

\begin{enumerate}

\item If $h^0(\Te_X(\a)),h^0(\Te_X(\a+\e))\geq 2$, we can apply the Clifford Theorem to the divisors $\a,\a+\e$. We obtain:
$$
Cliff(C)\geq \pi-2g+1-2-2+2=\pi-2g-1.
$$

\item If $h^0(\a)\leq 1$ then $h^0(\Te_X(\a+\e))\geq 2$. The Clifford Theorem for the divisor $\a+\e$ gives:
$$
Cliff(C)\geq deg(\a)+\pi-2g+1-2\geq \pi-2g-1.
$$

\item If $h^0(\a+\e)=1$, then $h^0(\Te_X(\a))\geq 2$. Now, the Clifford Theorem applies to $\a$:
$$
Cliff(C)\geq deg(\a)-e+\pi-2g+1-2\geq \pi-2g-1.
$$

\item If $h^0(\aa+\e)=0$, then $h^0(\Te_X(\a))\geq 3$. We know that:
$$
deg(\a)-h^0(\Te_X(\aa))\geq gon(X)-2.
$$
From this and the Lemma \ref{bound0}:
$$
Cliff(C)\geq 2gon(X)-4-e+\pi-2g+1+2\geq 2gon(X)-2.
$$ \qed

\end{enumerate}

\begin{cor}

Let $C$ be a smooth curve of genus $\pi$ with a double covering $\gamma:C\lrw X$ over a smooth curve of $C$ of genus $g$. If $\pi\geq 2g-1+2gon(X)$, then:
$$
gon(C)=2gon(X)\qquad \hbox{ and }\qquad Cliff(X)=2gon(X)-2.
$$
\end{cor}
{\bf Proof:} It is a consequence of the Theorem \ref{gonality}, \ref{clifford} and the Lemma \ref{bound1}. \qed

This result says nothing when $\pi< 2g-1+2gon(C)$. However, as an application of the Theorem \ref{lowdegree} we can construct double coverings of a generic curve $X$ of genus $g$ with prescribed gonality:

\begin{cor}

Let $X$ be a  generic curve of genus $g$. Then for any integer $d$ with $gon(X)\leq d\leq 2gon(X)$ there exists a double covering $\gamma:C\lrw X$ with $C$ a smooth curve such that:
$$
gon(C)=d, \qquad Cliff(C)=d-2.
$$

\end{cor}
{\bf Proof:} Since $X$ is generic, there is a base-point-free pencil of degree $d\geq gon(x)$. Now, it is sufficient to apply the Theorems \ref{lowdegree}, \ref{gonality} and \ref{clifford}. \qed


\begin{thebibliography}{77}

\bibitem{BaKe95}{\sc Ballico, E.; Keem, C.}
{\it On multiple coverings of irrational curves.}
Arch. Math {\bf N.65}, 151-160 (1995).

\bibitem{BaKePa04}{\sc Ballico, E.; Keem, C.; Park, S.}
{\it Double covering of curves.}
Proc. Amer. Math. Soc. 132, {\bf N. 11}, 3153-3158 (2004).

\bibitem{FuPe05}\textsc{Fuentes Garc\'{\i}a, L., Pedreira, M.}:
{\it The projective theory of ruled surfaces.}
{Note Mat., {\bf 24}, 25-63 (2005).} 

\bibitem{FuPe05*}\textsc{Fuentes Garc\'{\i}a, L., Pedreira, M.}:
{\it Canonical geometrically ruled surfaces.}
Math. Nachr. \textbf{278}, 240-257 (2005).

\bibitem{Ha77}{\sc Hartshorne, R. }
{\it Algebraic Geometry. }
GTM, 52. Springer--Verlag, 1977.

\bibitem{Ka84}{\sc Kani, E.}
{\it On Castelnuovo's equivalence defect.}
J. Reine Angew. Math. {\bf 352}, 24-70 (1984).

\bibitem{KeOh06}{\sc Keem, C.; Ohbuchi }
{\it On the Castelnuovo-Severi inequality for a double covering. }
Preprint.


\end{thebibliography}
\end{document}